\documentclass[12pt]{article}
 \usepackage[margin=0.8in]{geometry} 
\usepackage{amsmath,amsthm,amssymb,amsfonts,stmaryrd}
\usepackage{enumitem}
\usepackage{graphicx}
\usepackage{tikz}
\usepackage{bbm}
\usepackage{multicol}
\usepackage{array}
\usepackage{hyperref}
\usepackage{xcolor}
\usepackage{ocgx}
\usepackage{mathabx}
\usepackage{subfig}
 
 \usepackage{tikz}
\usetikzlibrary{patterns}

\newcommand{\RR}{\mathbb{R}}
\newcommand{\PP}{\mathbb{P}}

\newcommand{\ZZ}{\mathbb{Z}}

\newcommand{\norm}[1]{\|#1\|} 

\theoremstyle{definition}

 \newtheorem{theorem}{Theorem}

\usepackage{biblatex}
\begin{document}
\title{Stochastic comparisons of Tracy-Widom $\beta$ distributions}
\author{Virginia Pedreira}
\date{}
\maketitle

\begin{abstract}
    We give a stochastic comparison and ordering of the  Tracy-Widom distribution with parameter $\beta$. In particular, we show that as $\beta$ grows, the Tracy-Widom random variables get smaller modulo a multiplicative coefficient.
\end{abstract}

\section{Introduction}

In a sequence of papers, J. Baik, P. Deift,  K. Johansson and E. Rains ( \cite{lengthoflongest1999}, \cite{johansson2000}, \cite{symmetrizedrandompermutations2001}, \cite{algebraicaspectsofincreasing2001} and \cite{asymptoticsofmonotone2001}) define several versions of a last passage percolation model on the $\ZZ^2$ lattice and take the rescaled limit of those models. They assign weights to each vertex of the lattice and they run random weighted walks on the square $[N, N]^2$. They study the paths with the largest weight and find that, in the rescaled limit, those paths converge to the Tracy-Widom distribution. They also apply certain symmetries to the lattice and obtain last passage paths that follow those symmetries and, rescaled, converge to Tracy-Widom distributions with different parameters (1, 2 or 4). In those papers, they imply that there is an interpolation of those last passage paths that implies an interpolation of the limits. In this paper, we show how we can stochastically compare and order the different Tracy-Widom distributions with general parameter $\beta$, proving that interpolation in the general case. The rest of this section will explain the motivation and state the result. In section \ref{theorysection} we will explain the background material and in section \ref{proofssection} we will prove the result. 

 In \cite{symmetrizedrandompermutations2001}, Baik and Rains obtained the asymptotic fluctuations of the models mentioned above that we will now define. To each site $(i,j)\in\ZZ^2$ we assign a random variable $w(i,j)$. The random variables at each site are independent and identically distributed. We will denote a general up/right path as $\pi:(i,j)\nearrow(k,l)$, indicating its initial and final position. The weight or length of each path is the sum of the weights of the sites it visits. The goal is to describe the asymptotic length of the longest up/right path. We will apply three symmetries $T_{\boxasterisk}:\ZZ^2\rightarrow\ZZ^2$ onto the lattice. The identity symmetry will be named $T_{\boxvoid}$, the symmetry along the $y=x$ diagonal will be called $T_{\boxslash}$ and the symmetry along the other diagonal, $y=-x$ will be called $T_{\boxbackslash}$.
 
 Then, the length of the longest up/right path on a square with side length $N$ can be described as $G^{T_{\boxasterisk}}(N)=\sup_{\pi:p \nearrow q}\sum_{(i,j)\in\pi}w(T_{\boxasterisk}(i,j))$. We always take a square of size $N$ but taking into account that the "diagonals" of the lattice have to coincide with the diagonals of the square, the initial and final points $p$ and $q$ might be different for each symmetry. However, the points $p$ and $q$ always represent the lower left and upper right points in the square (but we can not always take the square $[0,N]\times[0,N]$).
 
 In this context, Baik and Rains proved that for each $x\in\RR$, 
 \[
 \lim_{N\rightarrow\infty}\PP\left(\frac{G^{T_{\boxasterisk}}(N)-aN)}{bN^{1/3}}\leq x\right)=F_*(x), 
 \]
 where the constants $a$ and $b$ depend on the distribution of the weights $w(i,j)$ and  for each symmetry $\boxasterisk=\boxvoid, \boxslash$ or  $\boxbackslash$, the function $F_*(x)$  is the cumulative distribution function of the Tracy-Widom 2, 4 and 1 respectively, as originally defined in \cite{TracyWidomGUE} and \cite{TracyWidomGOEGSE} by Tracy and Widom. We will name the random variable associated to $F_*$ as $L^{\boxasterisk}$.
 
We can compare $L^{\boxvoid}$ and $L^{\boxslash}$ by defining the last passage model as above on the square $[0,N]\times[0,N]$ and in the case of $T_{\boxslash}$ we symmetrize the half plane above the diagonal onto the lower half plane. This coupling gives us a simple comparison of $L^{\boxvoid}$ and $L^{\boxslash}$: in the case of the $\boxslash$ symmetry, we are taking the maximum of the up/right paths that stay in the upper half triangle while in the case of $\boxvoid$ the maximum is taken on all the up/right paths from the lower left corner to the upper right corner of the square. Since the weights in the upper half triangle are the same in both models,  $G^{\boxvoid}$ is larger than $G^{\boxslash}$. Therefore, $L^{\boxvoid}\geq L^{\boxslash}$.

Similarly, we can compare $L^{\boxvoid}$ with $L^{\boxbackslash}$. In this case, we define the last passage model in the square $[-N,0]\times[0,N]$. The model is shift invariant in the lattice and this square will allow us to couple both random variables. The symmetry $T_{\boxbackslash}$ acts by copying the triangle below the $y=-x$ diagonal onto the upper triangle symmetrically. We can see that every up/right path from the lower left corner of the square to the upper right corner of the square in the symmetrized lattice consists of two symmetric paths: the path form the lower left corner to the diagonal is then repeated symmetrically in the upper triangle. Therefore, the weight longest path is exactly twice the weight of the longest path from the lower left corner to the diagonal. This path to the diagonal be larger than the path from the lower left corner to the center of the square because the center of the square is in the diagonal, assuming that $N$ is even. The weight of this new path is equal to $G^{\boxvoid}(N/2)$ so 
\[
\frac{1}{2}G^{\boxbackslash}(N)=G^{\boxvoid}\left(\frac{N}{2}\right).
\]
After substracting the mean, rescaling and taking the limit, we obtain that 
\[
L^{\boxbackslash}\geq2^{2/3}L^{\boxvoid}.
\]

As mentioned before, the random variables  $L_{\boxasterisk}$ are distributed according to the Tracy-Widom distributions as defined in \cite{TracyWidomGOEGSE} and \cite{TracyWidomGUE}. In fact, $L_{\boxasterisk}$ is the rescaled limit of the largest eigenvalue of a Gaussian random matrix. In   \cite{sao} , Ram\'irez, Rider and Vir\'ag, propose a tridiagonal random matrix that depends on a parameter $\beta$ and whose spectrum distribution, called the $\beta$-ensemble, coincides with the Gaussian Ensembles (GO/U/SE) in the cases where $\beta$ is 1, 2 or 4. In that sense, they generalize the Tracy-Widom with parameter $\beta$ by taking the rescaled limit of the largest eigenvalue. We call those random variables as $TW_{\beta}$. This new definition differs slightly from the original one for the cases where $\beta$ is 1, 2 or 4. An explanation on the way the scaling differs in the two definitions can be found in the work of Bloemendal and Vir\'ag \cite{TWbetadefinitionbloemendal}. This slight difference means that, $L_{\boxbackslash}=TW_1$, $L_{\boxvoid}=TW_2$ and $L_{\boxslash}=2^{2/3}TW_4$. 

From the coupling, we see a pattern in these stochastic comparisons:
\[
TW_1\geq2^{2/3} \quad \text{ and } \quad TW_2\geq2^{2/3}TW_4.
\]
We will prove that this generalizes for the Tracy-Widom $\beta$ random variables defined originally, by Ram\'irez, Rider and Vir\'ag in \cite{sao}.

The main result in this work is the following:
\begin{theorem}
\label{theorem}
Let $\beta'>\beta>0$ and $\alpha>0$, then $TW_{\beta}\geq\alpha TW_{\beta'}$ if and only if $\left(\frac{\beta'}{\beta}\right)^{1/3}\leq\alpha\leq\left(\frac{\beta'}{\beta}\right)^{2/3}$.
\end{theorem}  

\begin{figure}[h]
    \centering
   \includegraphics[scale=1]{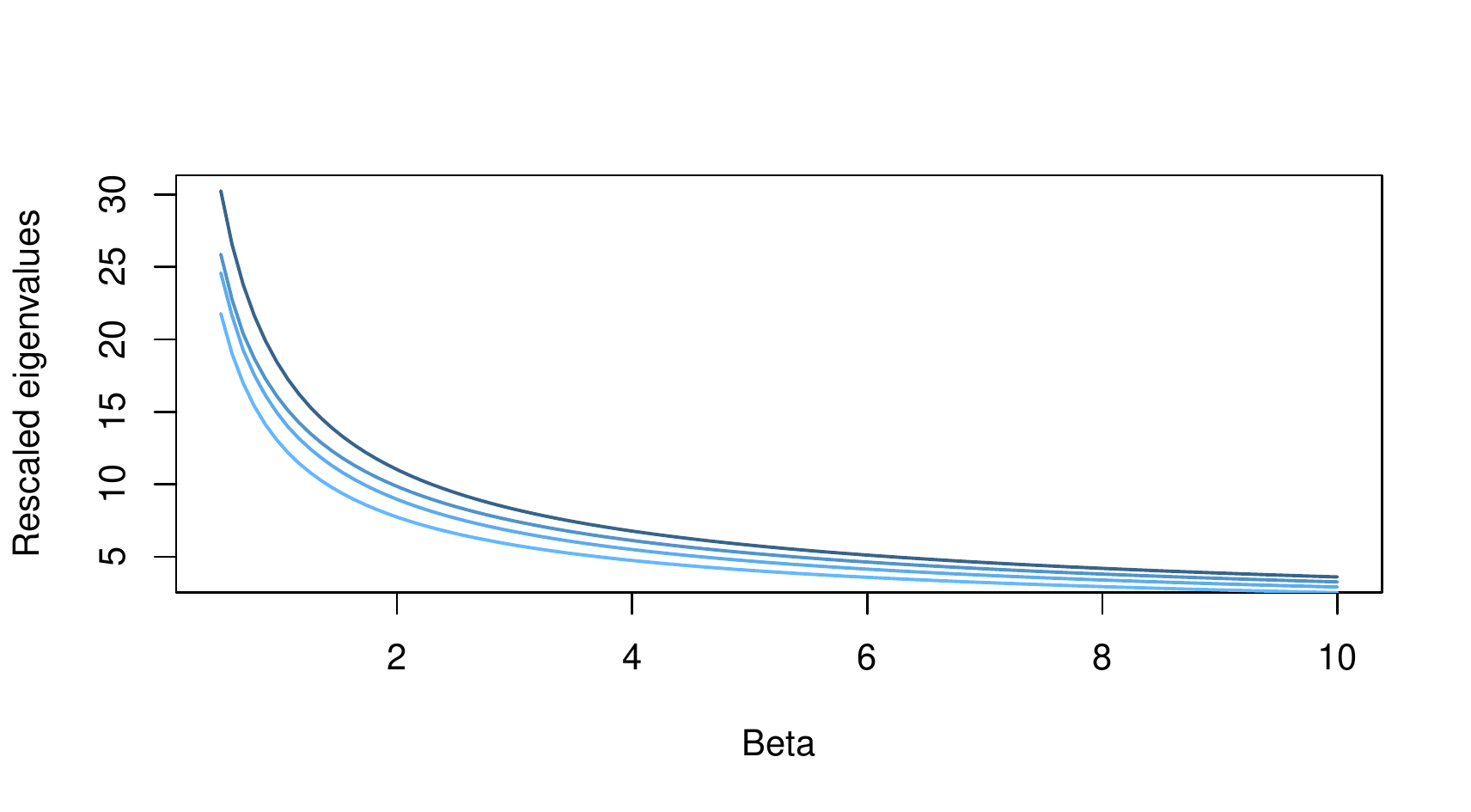}
    \caption{The first four rescaled eigenvalues $\frac{(\lambda_i(\beta)-2\sqrt{n})n^{1/6}}{\beta^{2/3}}$ of a 10x10 matrix distributed according to the $\beta$-ensemble plotted as functions of $\beta$. The colour gradient represents the order of the eigenvalues; the lighter the shade, the smaller the eigenvalue. The functions are decreasing. }
    \label{eigenvalues vs beta}
\end{figure}

\section{Stochastic Airy operator and Tracy Widom beta random variables}
\label{theorysection}
    For any $\beta>0$, consider the probability distribution 
\[
\PP_n^{\beta}(\lambda_1,\lambda_2,\dots,\lambda_n)=\frac{1}{Z_n^{\beta}}e^{-\beta\sum_{k=1}^n\lambda_k^2/4}\prod_{j<k}|\lambda_j-\lambda_k|^{\beta},
\]
where $\lambda_1\geq\lambda_2\geq\dots\geq\lambda_n$. When $\beta=1,2$ and 4, this distribution corresponds to the distribution of the joint density of eigenvalues of Gaussian orthogonal, unitary and symplectic ensembles respectively, or G(O/U/S)E, of random matrix theory. In \cite{sao}, Ram\'irez, Rider and Vir\'ag obtain the point process limit of the spectral edge of the general $\beta$-ensemble. In fact, the eigenvalues of the $\beta$-ensemble converge in distribution to the eigenvalues of a stochastic operator called the Stochastic Airy Operator (SAO):
\begin{equation}
\label{operator}
H_{\beta}=-\frac{d^2}{dx^2}+x+\frac{2}{\sqrt{\beta}}b_x'
\end{equation}
where $b'$ is the white noise. The operator is defined on the Hilbert space $L^*$, the space of continuous functions $f$ such that $f(0)=0$ and $\int_0^{\infty}(f'(x))^2+(1+x)f^2(x)dx<\infty$. The SAO acts on functions as a quadratic form is the following way: we decompose the Brownian motion in two terms $b=\overline{b}+(b-\overline{b})$ where $\overline{b}(x)$ is the average, 
\[
\overline{b}(x)=\int_x^{x+1}b_ydy.
\]
For every function $f\in L^*$, 
\[
\prec f,H_{\beta}f\succ=\int_0^{\infty}(f'(x))^2+xf^2(x)dx+\frac{2}{\sqrt{\beta}}\left(\int_0^{\infty}f^2(x)\overline{b}_x'dx-2\int_0^{\infty}f'(x)f(x)(\overline{b}_x-b_x)dx\right).
\]
The integrals are well defined and finite, see \cite{sao}. The definition looks more involved than it needs to be and that is because the last two integrals could be replaced, using integration by parts, with $\frac{-4}{\sqrt{\beta}}\int_0^{\infty}f(x)f'(x)b_xdx$ if this were a finite integral. If the function $f$ is compactly supported, then this simpler definition of the quadratic form works.

We use the variational characterization of the eigenvalues and eigenfunctions.
Then, the smallest eigenvalue, $\Lambda_0$, is defined as
\[
\Lambda_0=\inf\{\prec f,H_{\beta}f\succ : f\in L^*, \norm{f}_2=1\}.
\]
The infimum of the formula above is attained at an eigenfunction $f_0$ with corresponding eigenvalue $\Lambda_0$.   Functions of compact support are dense in $L^*$ (since functions in $L^*$ have the boundary condition $f(0)=0$) and the quadratic form $\prec\cdot,H_{\beta}\cdot\cdot\succ:(L^*)^2\rightarrow\RR$ is continuous as proved in \cite{sao} so we can take $\Lambda_0$ to be 
\[
\Lambda_0=\inf\{\prec f,H_{\beta}f\succ : f\in L^*, \norm{f}_2=1, f \text{ compactly supported}\}.
\]
The rest of the eigenvalues are defined recursively as
\[
\Lambda_k=\inf\{\prec f,H_{\beta}f\succ : f\in L^2, \norm{f}_2=1, f\perp f_0, \dots, f_{k-1}\}
\]
and the supremum is attained at an eigenfunction $f_k$. Functions of compact support are also $L^*$-dense in the orthogonal complement of an eigenspace. In fact, if a function $f$ in $L^*$ is orthogonal to $f_0,\dots,f_{k-1}$ we can choose a function $g$ that is $\varepsilon$ close to $f$ in $L^*$. Then, the function $\Tilde{g}=g-\sum_{i=0}^{k-1}\langle g,f_i\rangle f_i$ is a function in $L^*$  orthogonal to the eigenfunctions $f_i$. Since the function $f$ is also orthogonal to the eigenfunctions, we can rewrite $f$ as $f=f-\sum_{i=0}^{k-1}\langle g,f_i\rangle f_i$. Then, 
\[
\norm{f-\Tilde{g}}_{L^*}=\norm{f-g-\sum_{i=0}^{k-1}(\langle f,f_i\rangle-\langle g,f_i\rangle) f_i}_{L^*}\leq\norm{f-g}_{L^*}+\sum_{i=0}^{k-1}|\langle f,f_i\rangle-\langle g,f_i\rangle|\leq(k+1)\varepsilon
\]
using Cauchy-Schwarz inequality and the fact that $\norm{f}_2\leq\norm{f}_{L^*}$. 
The continuity of the quadratic form in $L^*$ means that we can restrict the definition to 
\[
\Lambda_k=\inf\{\prec f,H_{\beta}f\succ : f\in L^2, \norm{f}_2=1, f\perp f_0, \dots, f_{k-1},f \text{ compactly supported}\}.
\] 
 More details on this random operator and its eigenvalues can be found in \cite{sao}.

Then, (from \cite{sao}) $\Lambda_0\leq\Lambda_1\leq\dots\leq\Lambda_{k-1}$ are the $k$ lowest elements of the set of eigenvalues of the operator $H_{\beta}$ and the vector
\[
(n^{1/6}(2\sqrt{n}-\lambda_l))_{l=1,\dots,k}
\]
converges in distribution to $(\Lambda_0,\Lambda_1,\dots\Lambda_{k-1})$ as $n\rightarrow\infty$.

The rescaled limit of the largest eigenvalue of the $\beta$-ensembles mentioned earlier is distributed according the Tracy-Widom $\beta$, so we define the  Tracy-Widom $\beta$ distribution as the distribution of $-\Lambda_0$. In fact, 
\[
TW_{\beta}=-\Lambda_0
\]

There is a deterministic operator associated with the SAO which is the Airy Operator
\[
A=-\frac{d^2}{dx^2}+x.
\]
We can think of the Airy operator as the SAO with parameter $\beta=\infty$. 
    
\section{Proof of Theorem \ref{theorem}}
\label{proofssection}

The goal is to stochastically compare the eigenvalues of the SAO$_{\beta}$. Recall that the parameter $\beta$ only appears in the operator \ref{operator} as part of the coefficient of the random term, so the coupling used to obtain the comparison will consist of keeping the same source of randomness for all $\beta$.

There is a natural partial order on the space of self adjoint operators, the Loewner order: we say that two operators $A$ and $B$ are ordered $A\geq B$ if the operator $A-B$ is positive definite. We would like to establish an order on $\{H_{\beta}\}_{\beta\geq1}$. 

Assume that $\beta'>\beta>0$. Then, we will show that there exists a constant $c$ such that $cH_{\beta'}\geq H_{\beta}$. 



We will rescale the SAO equation as follows: take $y=sx$ where $s$ is a fixed positive number. We call $f_s(y):=f(y/s)=f(x)$. Then, doing the corresponding change of variables, we get
\begin{align*}
    \int_0^{\infty}\left(f'(x)\right)^2dx&=s\int_0^{\infty}\left(f'_s(y)\right)^2dy\\
    \int_0^{\infty}xf^2(x)dx&=\frac{1}{s^2}\int_0^{\infty}yf_s^2(y)dy\\
  \int_0^{\infty}f(x)f'(x)b_xdx&=\frac{1}{\sqrt{s}}\int_0^{\infty}f_s(y)f_s'(y)b_ydy
\end{align*}
using the Brownian scaling. Therefore, we have the following identity in distribution 
\[
H_{\beta}=-\partial_x^2+x+\frac{2}{\sqrt{\beta}}b_x'=-s\partial_y^2+\frac{1}{s^2}y+\frac{2}{\sqrt{s\beta}}b_y'=:H_{\beta}^s.
\]
 Let $\Lambda_k^{\beta,s}$ be the $k$-th eigenvalue of $H_{\beta}^s$. Then, $\Lambda_k^{\beta,s}=\frac{1}{s}\Lambda_k^{\beta}$

Let $\beta'>\beta>0$ and $\gamma=\sqrt{\frac{\beta'}{s\beta}}$. Then,
\[
\gamma H_{\beta'}=-\sqrt{\frac{\beta'}{s\beta}}\partial_y^2+\sqrt{\frac{\beta'}{s\beta}}y+\frac{2}{\sqrt{s\beta}}b_y'
\]
so 
\[
\gamma H_{\beta'}-H_{\beta}^s=-\left(\sqrt{\frac{\beta'}{s\beta}}-s\right)\partial_y^2+\left(\sqrt{\frac{\beta'}{s\beta}}-\frac{1}{s^2}\right)y.
\]

Notice that the Airy operator is positive definite since 
\[
\left\langle Af,f\right\rangle=\int_0^{\infty}-f''(x)f(x)+xf^2(x)dx=\int_0^{\infty}f'^2(x)+xf^2(x)dx\geq0.
\]
In fact, if we take the deterministic operator $A_{a,b}=-a\partial_x^2+bx$, we know that $A_{a,b}$ is positive definite if and only if both $a$ and $b$ are positive. Then, we need that 
\begin{align*}
   \sqrt{\frac{\beta'}{s\beta}}&\geq s\\
   \sqrt{\frac{\beta'}{s\beta}}&\geq\frac{1}{s^2}
\end{align*}
which happens if and only if 
\[
\sqrt{\frac{\beta}{\beta'}}\leq s^{3/2}\leq\sqrt{\frac{\beta'}{\beta}}.
\]
Notice that since $\beta'>\beta$, the inequality makes sense. In other words, we need 
\[
\left(\frac{\beta}{\beta'}\right)^{1/3}\leq s\leq\left(\frac{\beta'}{\beta}\right)^{1/3}.
\]
We conclude that $\gamma H_{\beta'}\geq H_\beta^s$ for all the coefficients $s$ in that range. Since the positive definite partial order implies an ordering of the eigenvalues, we have that $\gamma\Lambda_k^{\beta'}\geq\Lambda_k^{\beta,s}=\frac{1}{s}\Lambda_k$  for the same range of $s$. Let $\alpha=s\gamma$. The restriction on $\left(\frac{\beta}{\beta'}\right)^{1/3}\leq s\leq\left(\frac{\beta'}{\beta}\right)^{1/3}$ is equivalent to the restriction on $\alpha$ given by
\[
\left(\frac{\beta'}{\beta}\right)^{1/3}\leq \alpha=\sqrt{s}\sqrt{\frac{\beta'}{\beta}}\leq\left(\frac{\beta'}{\beta}\right)^{2/3}
\]

 We have proved that if $\left(\frac{\beta'}{\beta}\right)^{1/3}\leq\alpha\leq\left(\frac{\beta'}{\beta}\right)^{2/3}$, then $TW_{\beta}\geq\alpha TW_{\beta'}$. (Here, the inequality reverses because the $-\Lambda_0$ is distributed according to $TW_{\beta}$)

Notice that this proof gives a comparison of the whole spectrum of the Stochastic Airy Process and not only on the smallest eigenvalue which is distributed according to the Tracy-Widom distribution. In fact, if, as before, $\Lambda_0^{\beta}, \Lambda_1^{\beta}, \dots$ are the eigenvalues of $H_{\beta}$ in increasing order, 
\[
\alpha\Lambda_k^{\beta'}\geq\Lambda_k^{\beta},
\]
for any $k$, given that $\left(\frac{\beta'}{\beta}\right)^{1/3}\leq\alpha\leq\left(\frac{\beta'}{\beta}\right)^{2/3}$.

In the opposite direction, we can look at the tails of the $TW_{\beta}$ distribution and get from there a possible range of $\alpha$s. From Virág, Ramírez, Rider, \cite{sao} we get 
\begin{align*}
    \PP(TW_{\beta}>a)&=\exp\left(-\frac{2}{3}\beta a^{3/2}(1+o(1))\right)\\
    \PP(TW_{\beta}<-a)&=\exp\left(-\frac{1}{24}\beta a^3(1+o(1))\right).
\end{align*}
If $\alpha TW_{\beta'}\leq TW_{\beta}$, then $\PP(\alpha TW_{\beta'}>a)\leq\PP(TW_{\beta}>a)$ which means that 
\[
\exp\left(-\frac{2}{3}\beta'\frac{a^{3/2}}{\alpha^{3/2}}(1+o(1))\right)\leq\exp\left(-\frac{2}{3}\beta a^{3/2}(1+o(1))\right)
\]
or 
\[
\exp\left(-\frac{2}{3}a^{3/2}(\beta'/\alpha^{3/2}-\beta)(1+o(1))\right)\leq 1
\]
so $\beta'/\alpha^{3/2}-\beta\geq 0$ or equivalently, $\alpha\leq\left(\frac{\beta'}{\beta}\right)^{2/3}$. Doing a similar calculation with the left-hand tail, gives us that $\left(\frac{\beta'}{\beta}\right)^{1/3}\leq\alpha$ which is the same range that we found through the other method. This concludes the proof.

\printbibliography
\end{document}